\documentclass[reqno,letterpaper]{amsart}

\usepackage{amscd,amssymb} \usepackage[cmtip,arrow,curve,matrix]{xy}

\newcommand{\splitpageyesno}[2]{#2} \splitpageyesno{
\usepackage{splitpage} } { \usepackage{fullpage} }

 \numberwithin{equation}{section}
\newtheorem{Theorem}[equation]{Theorem}

\newtheorem{Lemma}[equation]{Lemma}
\newtheorem{Proposition}[equation]{Proposition}
\newtheorem{Corollary}[equation]{Corollary}

\theoremstyle{definition}

\newtheorem{Definition}[equation]{Definition}
\newtheorem{Remark}[equation]{Remark}
\newtheorem{Example}[equation]{Example}

\DeclareMathOperator{\Aut}{Aut}

\newcommand{\C}{\mathbb{C}} 
\newcommand{\cpinfty}{\C \! P^{\infty}}
\newcommand{\cp}{P}

\DeclareMathOperator*{\colim}{colim}

\newcommand{\eqdef}{\overset{\text{def}}{=}}

\newcommand{\Nu}{\mathcal{V}}

\renewcommand{\O}{\mathcal{O}}

\newcommand{\psb}[1]{[ \! [#1] \! ]}

\newcommand{\Qp}{\mathbb{Q}_{p}}

\newcommand{\R}{\mathbb{R}}

\newcommand{\sinf}{\Sigma^{\infty}} \DeclareMathOperator{\spec}{spec}
 \DeclareMathOperator{\spf}{spf}
 
\newcommand{\slot}{\,-\,}

\newcommand{\tensor}[1]{\underset{#1}{\otimes}}

\newcommand{\xra}[1]{\xrightarrow{#1}}

\newcommand{\Z}{\mathbb{Z}} 
\newcommand{\Zp}{\Z_{p}}

\newcommand{\antvg}{\ellg^{an}}
\newcommand{\bara}{\bar{a}}
\newcommand{\barm}{\bar{m}} 
\newcommand{\borel}[1]{#1_{G}}
\newcommand{\busix}{\BU{6}} 
\newcommand{\BU}[1]{BU\langle #1 \rangle}

\newcommand{\cochars}{\check{T}}
\newcommand{\cx}[1]{#1^{\C}}

\newcommand{\dirac}{D}

\DeclareMathOperator{\ellg}{Ell}
\newcommand{\einfty}{E_{\infty}} 
\newcommand{\Ell}[1]{\mathbf{#1}}

\DeclareMathOperator{\Ext}{Ext}

\newcommand{\fmlgpof}[1]{\widehat{#1}} 
\newcommand{\fs}[1]{+_{#1}}

\newcommand{\Gah}{\widehat{\mathbb{G}}_{a}}
\newcommand{\GpOf}[1]{P_{#1}}

\DeclareMathOperator{\Hom}{Hom} 

\newcommand{\hkr}{\Xi}

\newcommand{\Ihat}{\hat{I}} 
\DeclareMathOperator{\ind}{ind}

\newcommand{\KTate}{K_{\mathrm{Tate}}}

\newcommand{\mspin}{MSpin} 
\newcommand{\MU}[1]{MU\langle #1 \rangle}
\newcommand{\musix}{\MU{6}} 
\newcommand{\moeight}{MO\langle 8 \rangle}

\newcommand{\orbsigma}{\sigma_{\mathrm{orb}}}
\newcommand{\orbtvg}{\ellg_{orb}}

\newcommand{\point}[1]{#1_{+}} 
\newcommand{\pontthom}[1]{\tau (#1)}

\newcommand{\restr}[1]{|_{#1}} 
\newcommand{\red}[1]{r #1}
\DeclareMathOperator{\rank}{rank}

\newcommand{\SigmaClass}{U} 
\newcommand{\sigmagenus}{\sigma}

\newcommand{\sr}[1]{w (#1)}

\newcommand{\uln}[1]{\underline{#1}} 
\newcommand{\T}{\mathbb{T}}
\newcommand{\tTate}{t} 
\newcommand{\Tate}{\mathrm{Tate}}

\begin{document}

\title[Discrete torsion and the sigma genus] {Discrete torsion for the
supersingular orbifold sigma genus}

\author[Ando] {Matthew Ando} \email{mando@math.uiuc.edu}

\author[French] {Christopher~P.~French} \email{cpfrench@math.uiuc.edu}

\address{Department of Mathematics \\
The University of Illinois at Urbana-Champaign \\
Urbana IL 61801 \\
USA}

\thanks{We thank the participants in the BCDE seminar at UIUC,
particularly David Berenstein and Eric Sharpe, for teaching us the
physics which led to this paper.  Ando was supported by NSF grant
number DMS---0071482.  Part of this work was carried out while Ando
was a visitor at the Isaac Newton Institute for Mathematical Sciences.  He
thanks the Newton Institute for its hospitality.}

\date{Version 1.4, August 2003.}

\begin{abstract}
The first purpose of this paper is to examine the relationship between
equivariant elliptic genera and orbifold elliptic genera.  We apply
the character theory of \cite{HKR:ggc} to the Borel-equivariant genus
associated to the sigma orientation of \cite{AHS:ESWGTC} to 
define an orbifold genus for certain total quotient orbifolds and
supersingular elliptic curves.  We show that our orbifold genus is
given by the same sort of formula as the orbifold ``two-variable'' genus of
\cite{DMVV:egspsqs} and \cite{BL:egsvoeg}.  In the case of a finite
cyclic orbifold group,  we use the characteristic series for the
two-variable genus in the formulae of \cite{Ando:AESO} to define an
analytic equivariant genus in Grojnowski's equivariant elliptic
cohomology, and we show that this gives \emph{precisely} the orbifold
two-variable genus. The second purpose of this paper is to study the
effect of varying the $\busix$-structure in the Borel-equivariant
sigma orientation.  We show that varying the $\busix$ structure by a
class in $H^{3}(BG;\Z)$, where $G$ is the orbifold group, produces
discrete torsion in the sense of \cite{vafa:dt}.  This result was
first obtained by Sharpe \cite{Sharpe:RDDT}, for a different orbifold
genus and using different methods.
\end{abstract}

\maketitle

\section{Introduction}\label{sec:introduction}

Let $E$ be an even periodic, homotopy commutative ring spectrum, let
$C$ be an elliptic curve over $S_{E} = \spec \pi_{0}E$, and let $t$ be
an isomorphism of formal groups
\[
   t: \fmlgpof{C} \cong \spf E^{0} (\cpinfty),
\]
so that $\Ell{C} = (E,C,t)$ is an elliptic spectrum in the sense of
\cite{ho:icm,AHS:ESWGTC}.  In \cite{AHS:ESWGTC}, Hopkins, Strickland,
and the first author construct a map of homotopy commutative ring spectra
\[
   \sigma (\Ell{C}): \musix \xra{} E
\]
called the \emph{sigma orientation}; it is conjectured in
\cite{ho:icm} that this map is the restriction to $\musix$ of a
similar map $\moeight\to E.$

The sigma orientation is natural in the elliptic spectrum, and, if
$\KTate = (K\psb{q}, \Tate,\tTate)$ is the elliptic spectrum
associated to the Tate elliptic curve, then the map of homotopy rings
\begin{equation} \label{eq:1}
    \pi_{*} \musix \rightarrow \pi_{*}\KTate
\end{equation}
is the restriction from $\pi_{*}\mspin$ of the \emph{Witten genus}.
Explicitly, let $M$ be a Riemannian spin manifold, and let $\dirac$ be
its Dirac operator.  Let $T$ denote the tangent bundle of $M$.  If $V$
is a (real or complex) vector bundle over $M$, let $\cx{V}$ be the
complex vector bundle
\[
\cx{V} = V\tensor{\R}\C.
\]
If $V$ is a complex vector bundle, let $\red{V}$ be the reduced bundle
\[
    \red{V} = V - \uln{\rank V},
\]
and let
\[
   S_{t}V = \sum_{k\geq 0} t^{k} S^{k} V
\]
be the indicated formal power series in the symmetric powers of $V$.
The operation $V\mapsto S_{t}V$ extends to an exponential operation
\[
    K (X) \rightarrow K (X)\psb{t}
\]
because of the formula
\[
  S_{t} (V\oplus W) = (S_{t} V ) (S_{t} W).
\]
The \emph{Witten genus} of $M$ is given by the formula
\begin{equation} \label{eq:2}
w (M) = \ind \left(\dirac \otimes \bigotimes_{k\geq 1} S_{q^{k}}
(\red{\cx{T}})\right) \in \Z\psb{q},
\end{equation}
and the diagram
\[
\xymatrix{ {\pi_{*}\musix} \ar[r] \ar[dr]_-{\pi_{*}\sigma (\KTate )} &
{\pi_{*}\mspin}
 \ar[d]^{w} \\
& {\Z\psb{q}} }
\]
commutes \cite{AHS:ESWGTC}.

The Witten genus first arose in \cite{Witten:EllQFT}, where Witten
showed that various elliptic genera of a manifold $M$ are essentially
one-loop amplitudes of quantum field theories of closed strings moving
in $M$.  Locally on $M$ the quantum field theory associated to $w (M)$
is a conformal field theory, and the obstruction to assembling a
conformal field theory globally on $M$ is $c_{2}M.$ \footnote{There
are various ways to understand this obstruction
(\cite{Witten:EllQFT,bm:gd4ccI,GMS:Gerbes,Ando:AESO}).}  This gives a
physical proof that, if $c_{2}M = 0$, then $w (M)$ is the
$q$-expansion of a modular form.

Suppose that $M$ is an $SU$-manifold.  The formula~\eqref{eq:2} shows
that the Witten genus is an invariant of the spin structure of $M$.
On the other hand the sigma orientation depends on a choice of
$\busix$ structure, that is, a lift in the diagram
\[
\xymatrix{ & {\busix} \ar[d]
\\
{M} \ar@{-->}[ur] \ar[r] & {BSU} \ar[r]^{c_{2}} & {K (\Z,4).}  }
\]
It is an interesting problem to understand how the orientation depends
on this choice.  The fibration sequence
\[
      K (\Z,3) \xrightarrow{\iota} \busix \rightarrow BSU \xra{c_{2}}
K (\Z,4)
\]
shows that a lift exists precisely when $c_{2} (M) = 0$, and that the set
of lifts is a quotient of $H^{3} (M;\Z)$.

The dependence on the choice of lift appears to have an explanation in
string theory.  The action for the QFT described by the Witten genus
is a function on the space of maps
\[
   X: \Sigma \to M
\]
of $2$-dimensional surfaces $\Sigma$ to $M$.  If the theory is
anomaly-free, that is, if $c_{2}M = 0$, then one is free to add to the
action a term of the form
\[
   \int_{\Sigma} X^{*}B,
\]
where $B\in \Omega^{2}M$ is a differential $2$-form on $M$ (called the
``B-field''), provided that
\[
    H = dB
\]
is an \emph{integral} three-form.  It seems clear that the physics of
the $B$-field should account for the variation in the sigma
orientation from at least torsion classes in $H^{3} (M,\Z)$.

In this paper we provide some evidence for this assertion.  We show
that varying the $\busix$ structure in the \emph{orbifold} sigma genus
of a supersingular elliptic curve produces the phenomenon known as
\emph{discrete torsion}, so named by Vafa in \cite{vafa:dt}.  Eric
Sharpe has shown that discrete torsion arises from the action of the
orbifold group on the $B$-field.  Putting Sharpe's results together
with ours suggests that, indeed, the $B$-field is the physical
reflection of the choice of $\busix$ structure.

More precisely, suppose that $M$ 
is a complex manifold with an action by a group $G$, and suppose that $V$ is a complex $G$-vector bundle
over $M$.  Let $T$ denote the tangent bundle of $M$.  If $X$ is a
space, let $X_{G}$ denote the Borel construction $EG\times_{G} X$.  If
\begin{align*}
      c_{1} (T_{G})  & = c_{1} (V_{G}) \\
      c_{2} (T_{G}) & = c_{2} (V_{G}),
\end{align*}
then there is a lift in the diagram
\begin{equation} \label{eq:18}
\xymatrix{ & {\busix}
 \ar[d] \\
{M_{G}} \ar[r]_{V_{G} - T_{G}} \ar@{-->}[ur]^{\ell} & {BSU,} }
\end{equation}
and a choice of lift gives a Thom class
\[
\SigmaClass  (M,\ell,\Ell{C})_{G} \in E^{0} (V_{G} - T_{G}).
\]
The relative zero section together with the Pontrjagin-Thom
construction provide a map
\[
\pontthom{V}_{G}: E^{0} (V_{G} - T_{G}) \rightarrow E^{0} (-T_{G}) \rightarrow
E^{0}(BG),
\]
and
\[
   \pontthom{V}_{G} (\SigmaClass (M,\ell,\Ell{C})_{G}) \in E^{0} (BG)
\]
is the (Borel) equivariant sigma genus of $M$ twisted by $V$ (see
\S\ref{sec:borel-equiv-sigma}).

To get from it an ``orbifold'' genus taking its values in $E^{0}$, we
use the character theory of
Hopkins, Kuhn, and Ravenel (\cite{HKR:ggc}; see also
\S\ref{sec:character-theory}).
It associates to a pair $(g,h)$ of commuting elements of $G$ a ring
homomorphism
\[
    \hkr_{g,h}: E^{0} (BG) \rightarrow D,
\]
where $D$ is a complete local $E$-algebra which depends on the formal
group of the spectrum $E$.  It turns out that the quantity
\begin{equation} \label{eq:17}
 \orbsigma (M,\ell,\Ell{C})_{G} =
\sum_{gh = hg} \hkr_{g,h} \pontthom{V}_{G} (\SigmaClass  (M,\ell,\Ell{C})_{G})
\end{equation}
takes its values in $E^{0}$; we call it the \emph{orbifold sigma
genus} of $M$ twisted by $V$ (see \S\ref{sec:orbifold-sigma-genus}).

There is already an extensive literature on the subject of ``orbifold
elliptic genera'', particularly the ``two-variable'' elliptic genus of
\cite{MR91e:57059,EOTY}; see for example \cite{DMVV:egspsqs,BL:egsvoeg}.  
In \S\ref{sec:orbifold-sigma-genus}, we show that the formula \eqref{eq:17}
is formally analogous to the formula for the orbifold two-variable
genus.   It is difficult to make a more precise comparison between our two
situations, because we work with the Borel-equivariant elliptic
cohomology associated to a supersingular elliptic curve, which is a
highly completed situation.  

In order to locate the orbifold two-variable genus more precisely in
the setting of equivariant elliptic cohomology, we consider in
\S\ref{sec:comp-with-analyt} the case of a finite 
cyclic group $G = \T[n] \subset \T$.  We use the principle suggested
by Shapiro's Lemma to define
\[
        E_{G} (X) \eqdef E_{\T} (\T\times_{G}X),
\]
where $E_{\T}$ is the uncompleted analytic equivariant elliptic
cohomology of Grojnowski.  We adapt the formulae in \cite{Ando:AESO},
which descends from \cite{Rosu:Rigidity,AndoBasterra:WGEEC}, to write
down an euler class in $E_{G} (X)$.  The associated genus $\antvg
(M,G)$ takes 
its value in $\Gamma (E_{G} (*))\cong \Gamma (\O_{C[n]})\cong
\C^{G\times G}$, and we prove the following.   

\begin{Theorem}
Summing the analytic equivariant two-variable genus over the torsion
points of the elliptic curve gives the orbifold two-variable genus:
more precisely, we have
\[
   \orbtvg (X,G)   =  \frac{1}{|G|}   \sum_{gh=hg} \antvg (M,G,g,h).
\]
\end{Theorem}

We were pleased to be able to
confirm that orbifold elliptic genera are so simply obtained from
equivariant elliptic genera.  It would be interesting to use this
observation to investigate more subtle properties of orbifold genera,
such as, for example, the ``McKay correspondence'' of Borisov and
Libgober.

The rest of the paper is devoted to the study of the dependence of
the orbifold sigma genus on the choice $\ell$ of $\busix$ structure
in \eqref{eq:18}.  Suppose that we have chosen an element $u \in H^{3} (BG;\Z)$,
represented as a map
\[
    u: BG \rightarrow K (\Z,3).
\]
If $\pi: M_{G} \rightarrow BG$ denotes the projection in the Borel
construction, then we obtain an element
\[
      \pi^{*}u = u \pi \in H^{3} (M_{G};\Z),
\]
and $\ell + \iota u \pi$ is another $\busix$-structure on $V_{G} -
T_{G}$.

In \S\ref{sec:cocyle}, we use the character theory and the sigma
orientation to associate to $u$ an alternating bilinear map
\[
 \delta = \delta (u,\Ell{C}): G_{2} \to D^{\times},
\]
where $G_{2}$ denotes the set of pairs of commuting elements of $G$ of
$p$-power order. In \S\ref{sec:formula} we obtain the

\begin{Theorem} \label{t-th-main-formula-intro}
The orbifold sigma genus associated to the $\busix$ structure
$\ell+\iota u \pi$ is related to the equivariant sigma genus
associated to $\ell$ by the formula
\[
         \orbsigma (M,\ell+ \iota u \pi,\Ell{C})_{G} = \sum_{gh = hg}
\delta(g,h) \hkr_{g,h} \pontthom{V}_{G} (\SigmaClass (M,\ell,\Ell{C})_{G}).
\]
\end{Theorem}

In \cite{vafa:dt}, Vafa observed that if
\[
       \phi = \sum_{gh=hg} \phi_{g,h}
\]
is an orbifold elliptic genus associated to a theory of strings on
$M$, and if $c=c (g,h)$ is a $2$-cocycle with values in $U (1)$, then
\begin{equation}\label{eq:11}
      \sum_{gh=hg} c (g,h)\phi_{g,h}
\end{equation}
is again modular; he called this phenomenon ``discrete torsion''.
Eric Sharpe \cite{Sharpe:RDDT} showed that the genus \eqref{eq:11}
arises from adding a $B$-field
\[
    B\in \Omega^{2} (M_{G})
\]
such that
\[
    dB = [c] \in H^{3} (M_{G};\Z),
\]
where $[c]$ is the cohomology class in $M_{G}$ obtained from $c$ by
pulling back along $M_{G}\to BG.$

Our result shows that varying the $\busix$-structure of $M_{G}$ by an
element $u \in H^{3} (BG)$ has a similar effect on the orbifold sigma
genus.  When $G$ is an abelian of order dividing $n=p^{s}$, the map
$\delta$ may be viewed as a two-cocycle on $G$ with 
values in $D^{\times}[n]\cong \Z/n$, and as such it represents a
cohomology class in $H^{2} (BG;\Z/n)\cong H^{3} (BG)$.  It is not
quite the cohomology class $u$: instead, as we shall see in
\S\ref{sec-non-abelian}, if $c$ is a $2$-cocycle representing $u\in
H^{3} (BG;\Z)\cong H^{2} (BG;\Z/n)$, then
\[
       \delta(g,h) = c (g,h) - c (h,g).
\]

\section{The sigma orientation and the sigma genus}
\label{sec:sigma-orient-sigma}

In this section we recall some results from \cite{AHS:ESWGTC}.

\begin{Definition}\label{def-elliptic-spectrum}
An \emph{elliptic spectrum} consists of \begin{enumerate} \item an
even, periodic, homotopy commutative ring spectrum $E$ with formal
group $P_E=\spf E^{0}\cpinfty$ over $\pi_0 E$; \item a generalized elliptic
curve $C$ over $\pi_{0}E$; \item an isomorphism $t:P_E\xra{}\fmlgpof{C}$
of $P_E$ with the formal completion of $C$.  \end{enumerate} A
\emph{map} $(f,s)$ of elliptic spectra $\mathbf{E_{1}} =
(E_1,C_1,t_1)\xra{} \mathbf{E_{2}}= (E_2,C_2,t_2)$ consists of a map
$f:E_1\xra{}E_2$ of multiplicative cohomology theories, together with
an isomorphism of elliptic curves
\[ C_{2} \xra{s} (\pi_0 f)_{\ast} C_{1}, \]
extending the induced isomorphism of formal groups.
\end{Definition}

\begin{Theorem} 
An elliptic spectrum $\Ell{C}= (E,C,t)$ determines a map
\[
      \sigma (\Ell{C}): \musix \rightarrow E
\]
of homotopy-commutative ring spectra.  The association $\Ell{C}\mapsto
\sigma (\Ell{C})$ is modular, in the sense that if
\[
(f,s): \Ell{C_{1}} \rightarrow \Ell{C_{2}}
\]
is a map of elliptic spectra, then the diagram
\[
\xymatrix{ {\musix} \ar[r]^-{\sigma (\Ell{C_{1}})} \ar[dr]_-{\sigma
(\Ell{C_{2}})} &
{E_{1}} \ar[d]^-{f}\\
& {E_{2}} }
\]
commutes up to homotopy.  If $\KTate = (K\psb{q},\Tate,\tTate)$ is the
elliptic spectrum associated to the Tate curve, then the diagram
\[
\xymatrix{ {\musix} \ar[r] \ar[dr]_{\sigma (\KTate)} & {MSpin}
\ar[d]^{w}
\\
& {K\psb{q}} }
\]
commutes, where $w$ is the orientation associated to the Witten genus.
\end{Theorem}

\section{The sigma genus} \label{sec:sigma-genus}

\begin{Definition} Let $W$ be a virtual complex vector bundle on a
space $M$.  A \emph{$\busix$-structure} on $W$ is a map
\[
   \ell: M \xra{} \busix
\]
such that the composition
\[
      M \xra{\ell} \busix \xra{} BU
\]
classifies $\red{W}$.
\end{Definition}

Now let $M$ be a connected compact closed manifold with complex tangent bundle
$T$, and let $V$ be another complex vector bundle on $M$.  Let 
\[
    d = 2\rank_{\C} T-V
\]
Let
\[
   \pontthom{V}: S^{0} \xra{P-T} M^{-T} \xra{\zeta} M^{V-T}
\]
be the composition of the Pontrjagin-Thom map with the
relative zero section.

If
\[
    \ell: M \to \busix
\]
is a $\busix$-structure on $V-T$, and if $\Ell{C} = (E,C,t)$ is an
elliptic spectrum, let $\SigmaClass (M,\ell, \Ell{C}) \in E^{-d}
(M^{V-T})$ be the class given by the map
\[
   \SigmaClass (M,\ell,\Ell{C}): \Sigma^{d} M^{V-T} \xra{\ell} \musix
\xra{\sigma (\Ell{C})} E.
\]

\begin{Definition} \label{def-sigma-genus}
The \emph{sigma genus} of $\ell$ in $\Ell{C}$ is the element
\[
   \sigmagenus (M,\ell,\Ell{C}) \eqdef \pontthom{V}^{*} (\SigmaClass
(M,\ell,\Ell{C})) \in E^{-d} (S^{0}) = \pi_{d} E.
\]
\end{Definition}

\begin{Example} \label{ex-1}
Suppose that $c_{1}T = 0 = c_{2}T$, so that $T$ itself admits a
$\busix$-structure, say $\ell: M\to \busix$, and $d=2\dim M$.  Then we
have a Thom isomorphism
\begin{align*}
      E^{0} (M) & \cong E^{-d} (M^{-T}) \\
             1 & \mapsto \SigmaClass (M,\ell, \Ell{C}).
\end{align*}
and the usual Umkehr map $\pi_{!}^{M}$ associated to the projection
\[
    \pi^{M}: M\to *
\]
is the composition
\[
\pi_{!}^{M}: E^{0} (M) \xra{\cong} E^{-d} (M^{-T}) \xra{P-T} E^{-d}
(S^{0}) \cong \pi_{d}E.
\]
Thus
\[
    \sigmagenus (M,\ell,\Ell{C}) = \pi_{!}^{M} (1) = \pi_{d} (\sigma
(\Ell{C})) ([M]) \in \pi_{d}E
\]
is just the genus of $M$ with $\busix$-structure $\ell$, associated to
the sigma orientation
\[
    \musix \xra{\sigma (\Ell{C})} E.
\]
\end{Example}

\section{The Borel-equivariant sigma genus}
\label{sec:borel-equiv-sigma}

Now suppose that $G$ is a compact Lie group, and, if $X$ is a space,
let $\borel{X}$ denote the Borel construction
\[
     \borel{X} \eqdef EG\times_{G} X.
\]
Suppose that $G$ acts on the compact connected manifold $M$, that $V$
is an equivariant complex vector bundle, and that
\[
    \ell: \borel{M} \to \busix
\]
is a $\busix$-structure on the bundle $\borel{V} - \borel{T}$.  Since
$\borel{T}$ is the bundle of tangents along the fiber of
\[
    \borel{M} \rightarrow BG,
\]
we have a Pontrjagin-Thom map
\[
   \point{BG}\xrightarrow{P-T} (\borel{M})^{-\borel{T}},
\]
and so a map
\[
\pontthom{V}_{G}: \point{BG}\xra{P-T} (\borel{M})^{-\borel{T}}
\xra{\zeta} (\borel{M})^{\borel{V} - \borel{T}}.
\]

Let $\SigmaClass (M,\ell,\Ell{C})_{G} \in E^{-d}
(\borel{M}^{\borel{V} - \borel{T}})$ be given by the map
\[
   \SigmaClass (M,\ell,\Ell{C})_{G}: \Sigma^{d}(\borel{M})^{\borel{V}
- \borel{T}} \xra{\ell} \musix \xra{\sigma (\Ell{C})} E.
\]
\begin{Definition} \label{def-equiv-sigma-genus}
The \emph{(Borel) equivariant sigma genus} of $\ell$ in $\Ell{C}$ is
the element
\[
   \sigmagenus (M,\ell,\Ell{C})_{G} \eqdef \pontthom{V}_{G} (\SigmaClass
(M,\ell,\Ell{C})_{G}) \in E^{-d} (BG).
\]
\end{Definition}

\section{Character theory} \label{sec:character-theory}

The equivariant sigma genus described in \S\ref{sec:borel-equiv-sigma}
is not so familiar, because $E^{*} (BG)$ is not.  In this section we
review the character theory of \cite{HKR:ggc}, which gives a sensible
way to understand $E^{*} (BG)$.  In the next section, we apply the
character theory to produce the orbifold sigma genus from the
equivariant sigma genus; as we shall see, it is given by the same sort
of formula as those for ``orbifold elliptic genera'' in for example
\cite{DMVV:egspsqs,BL:egsvoeg}

We suppose that $E$ is an even periodic ring spectrum, and that
$\pi_{0}E$ is a complete local ring of residue characteristic $p>0$.
We write $\cp$ for $\C P^{\infty}$, so $P_{E}=\spf E^{0}\cp$ is the
formal group of $E$.  We assume that $\GpOf{E}$ has finite height $h$.

Let $\Lambda_{\infty} = (\Zp)^{h}$, and for $n\geq 1$, let
$\Lambda_{n} = \Lambda_{\infty} / p^{n} \Lambda_{\infty}$.  If $A$ is
an abelian group, let $A^{*}\eqdef \hom (A,\C^{\times})$ denote its
group of complex characters, so for example $\Lambda_{\infty}^{*}
\cong (\Qp/\Zp)^{h}\cong (\Z[\tfrac{1}{p}]/\Z)^{h}$.

Each $\lambda\in \Lambda_{n}^{*}$ defines a map
\[
    B\Lambda_{n} \xrightarrow{B\lambda} \cp.
\]
Choose a coordinate $x\in E^{0}\cp.$ For each $\lambda \in
\Lambda_{n}^{*}$, let
\[
   x (\lambda) = (B\lambda)^{*}x \in E^{0}B\Lambda_{n}.
\]
Let $S\subset E^{0}B\Lambda_{n}$ be the multiplicative subset
generated by $\{x (\lambda) | \lambda \neq 0 \}$.  Let $L_{n} =
S^{-1}E^{0}B\Lambda_{n}$, and let $D_{n}$ be the image of
$E^{0}B\Lambda_{n}$ in $L_{n}.$ In other words, $D_{n}$ is the
quotient of $E^{0}B\Lambda_{n}$ by the ideal generated by annihilators
of euler classes of non-zero characters of $\Lambda_{n}$.  It is clear
that $L_{n}$ and $D_{n}$ are independent of the choice of coordinate
$x$.

Now suppose that $G$ is a finite group. Let
\[
      \alpha: \Lambda_{n} \to G
\]
be a homomorphism: specifying such $\alpha$ is equivalent to
specifying an $h$-tuple of commuting elements of $G$ of order dividing
$p^{n}$.

\begin{Definition} \label{def-char-map}
The \emph{character map associated to $\alpha$} is the ring
homomorphism
\[
     \hkr_{\alpha}: E^{0}BG \xra{E^{0}B\alpha} E^{0}B\Lambda_{n}
\xra{} D_{n}.
\]
\end{Definition}

One may check directly from the definition that the map $\Lambda_{n+1}
\twoheadrightarrow \Lambda_{n}$ induces maps
\begin{align*}
   D_{n} &\rightarrow  D_{n+1}\\
   L_{n} &\rightarrow L_{n+1}.
\end{align*}
Let
\begin{align}
    D = \colim_{n} D_{n} \notag \\
    L = \colim_{n} L_{n}. \label{eq:6}
\end{align}
Since $G$ is finite, any homomorphism
\[
    \alpha: \Lambda \rightarrow G
\]
factors as
\[
    \alpha: \Lambda \xra{} \Lambda_{n} \xra{\alpha_{n}} G
\]
for sufficiently large $n$, and we may unambiguously attach a
character homomorphism
\[
    \hkr_{\alpha}: E^{0}BG \xra{} D
\]
such that, for sufficiently large $n$, the diagram
\[
\xymatrix{ {E^{0}BG} \ar[r]^-{\hkr_{\alpha_{n}}}
\ar[dr]_-{\hkr_{\alpha}} & {D_{n}} \ar[d]
\\
& {D} }
\]
commutes.

A great deal is known about the ring $D_{n}$, because it turns out
\cite{AHS:hinfty} that $\spf D_{n}$ is the scheme of \emph{level
$\Lambda_{n}^{*}$-structures} on the $\GpOf{E}$.  
For example, it is
easy to check that the action of $\Aut (\Lambda_{n})$ on $E^{0}
B\Lambda_{n}$ induces an action of $\Aut (\Lambda_{n})$ on $D_{n}.$
Using the description of $D_{n}$ in terms of level structures, one may
prove the following.

\begin{Proposition}
The ring $D_{n}$ is finite and faithfully flat over $E$. If $\GpOf{E}$
is the universal deformation of a formal group of height $h$ (i.e. if
$E$ is a Morava $E$-theory), then $D_{n}$ is a complete Noetherian
local domain, and in that case, and in general if $p$ is regular in
$\pi_{0}E$, then $L_{n} = \frac{1}{p}D_{n}$.  The structural map
\[
   E^{0} \rightarrow D_{n}
\]
identifies $E^{0}$ with the $\Aut (\Lambda_{n})$-invariants in
$D_{n}$.
\end{Proposition}

\begin{proof}
\cite{Drinfeld:EM} or \cite{Strickland:FiniteSubgps}.
\end{proof}

If $w \in \Aut (\Lambda_{n})$ and $\alpha: \Lambda_{n}\to G$, then we
have two homomorphisms from $E^{0} (BG)$ to $D_{n}$, namely
$w\hkr_{\alpha}$ and $\hkr_{\alpha w}$.

\begin{Lemma} 
The diagram
\[
\xymatrix{ {E^{0}BG} \ar[r]^-{\hkr_{\alpha}} \ar[dr]_-{\hkr_{\alpha
w}} & {D_{n}} \ar[d]^{w}
\\
& {D_{n}} }
\]
commutes. \qed
\end{Lemma}


\begin{Corollary} \label{t-co-sum-is-in-E}
The expression
\[
     \sum_{\alpha:\Lambda_{n}\to G} \hkr_{\alpha}
\]
defines an additive map
\[
   E^{0}BG \rightarrow \pi_{0}E.
\] \qed
\end{Corollary}

In the case that the height of $\GpOf{E}$ is two, the sum is over all
pairs of commuting elements of $G$ of $p$-power order.  If $G$ is a
$p$-group, then we write
\[
     \sum_{gh=hg} \hkr_{g,h}: E^{0} BG \xra{} \pi_{0} E
\]
for the map in the Corollary.

\section{The orbifold sigma genus} \label{sec:orbifold-sigma-genus}

There has been much study of the orbifold version of the two-variable
elliptic genus of \cite{EOTY}; see for example
\cite{DMVV:egspsqs,BL:egsvoeg}.  In 
this section we introduce an orbifold version of the sigma genus, in
the case of a supersingular elliptic curve.  Our definition is
intentionally as simple as possible: we consider only total quotient
orbifolds, and then extract the orbifold sigma genus from the Borel
genus using the map of Corollary \ref{t-co-sum-is-in-E}.

Explicitly, suppose that $G$ is a finite group acting on a manifold
$M$ with complex tangent bundle $T$, that $V$ is an equivariant
complex vector bundle, and that
\[
   \ell: \borel{M} \to \busix
\]
is a $\busix$-structure on the bundle $\borel{V} - \borel{T}$.

Let $C$ be the universal deformation of a supersingular elliptic curve
over a perfect field of characteristic $p>0$, and let $\Ell{C} =
(E,C,t)$ be the associated elliptic spectrum.

\begin{Definition} \label{def-orbifold-genus}
The \emph{orbifold sigma genus} of $\ell$ in $\Ell{C}$ is the element
\begin{equation} \label{eq:5}
   \orbsigma (M,\ell,\Ell{C})_{G} \eqdef \sum_{gh=hg} \hkr_{g,h}
\sigmagenus (M,\ell,\Ell{C})_{G} \in \pi_{-d} E.
\end{equation}
\end{Definition}

The rest of this section is devoted to showing that the formula
\eqref{eq:5} is formally analogous to the formula for the orbifold
two-variable genus.  In section \ref{sec:comp-with-analyt}, we
show that, in the case of a finite cyclic group, the orbifold two-variable
genus is precisely the genus in 
Grojnowski's circle-equivariant elliptic cohomology obtained from the
characteristic series defining the two-variable genus by following the
construction of \cite{Rosu:Rigidity,AndoBasterra:WGEEC,Ando:AESO}.
These sections are logically independent of the
discussion of discrete torsion and the proof of Theorem
\ref{t-th-main-formula-intro}, and readers interested primarily in
that formula may prefer to skip to section~\ref{sec:cocyle}.

Our comparison in this section is based on the  
analogue of the formula \eqref{eq:5} in the case of a genus given by a
complex orientation
\[
     t: MU \xra{} E,
\]
so that $E$ has Thom classes and Umkehr maps for complex vector
bundles.

If $M$ is a compact manifold of real dimension $d$, then we write
$\pi^{M}$ for the projection
\[
   M\rightarrow *,
\]
and $\pi_{t}^{M}$ for the Umkehr map
\[
   \pi_{t}^{M} : E^{*} (M) \rightarrow \pi_{d-*} E.
\]
This Umkehr map is often denoted $\pi^{M}_{!}$; our notation
emphasizes the dependence on the orientation $t$.  In any case, the
genus associated to $t$ is the map
\[
    \Phi^{t}: \pi_{*} MU \rightarrow \pi_{*} E
\]
given by the formula
\[
   \Phi^{t} (M) = \pi_{t}^{M} (1).
\]

If a compact Lie group $G$ acts on $M$, then we write $\pi^{M,G}$ for
the projection
\[
     \pi^{M,G}: \borel{M} \rightarrow BG,
\]
and $\pi^{M,G}_{t}$ for the associated Umkehr map
\[
     \pi^{M,G}_{t}: E^{0} (\borel{M}) \xra{} E^{-d} (BG).
\]
If $G$ is a finite $p$-group, then the analogue of our
formula~\eqref{eq:5} is the quantity
\[
     \Phi^{t}_{\mathrm{orb}} (M) = \sum_{gh=hg} \hkr_{g,h}
\pi^{M,G}_{t} (1).
\]

If $g$ and $h$ are commuting elements of $G$, let
\[
      M^{(g,h)} \eqdef M^{g} \cap M^{h}
\]
be the subset of $M$ fixed by both $g$ and $h$.  Let $\Nu (g,h)$ be a
normal bundle of $M^{(g,h)}$ in $M$.

We view $(g,h)$ as a homomorphism
\[
   \Lambda_{n} \xra{(g,h)} G:
\]
this makes $\Lambda_{n}$ act trivially on $M^{(g,h)}$, and we let
\[
e_{t} (\Nu (g,h)) \in E^{*} (B\Lambda_{n})\otimes E^{*} (M^{(g,h)}) \cong
E^{*} (E\Lambda_{n}\times_{\Lambda_{n}} M^{(g,h)})
\]
be the $\Lambda_{n}$-equivariant euler class of $\Nu (g,h)$ in the
orientation $t$.

Recall from~\eqref{eq:6} that $L=\colim L_{n}$ is the colimit of the
rings $L_{n}$ obtained from $E^{*}B\Lambda_{n}$ by inverting the euler
classes of non-trivial characters of $\Lambda$.

\begin{Proposition}
Suppose that $(g,h) \neq (0,0)$.
\begin{enumerate}
\item The euler class $e_{t} (\Nu (g,h))$ is a unit of
\[
L\otimes_{E^{*}} E^{*} (M^{(g,h)}).
\]
\item The quantity
\[
     1 \otimes \pi_{t}^{M^{(g,h)}} \left( \frac{1}{e_{t} (\Nu (g,h))}
\right)
\]
lies in the subring $D\subset L$.
\item As elements of $D$ we have
\[
    \hkr_{g,h} (\pi_{t}^{M,G} (1)) = 1 \otimes \pi_{t}^{M^{(g,h)}}
\left( \frac{1}{e_{t} (\Nu (g,h))}\right).
\]
\end{enumerate}
Thus the orbifold genus of $M$ associated to $t$ is given by the
formula
\begin{equation} \label{eq:4}
          \Phi^{t}_{\mathrm{orb}} (M) = \Phi^{t} (M) +
\sum_{\substack{gh=hg \\ (g,h)\neq (0,0)}} 1\otimes
\pi_{t}^{M^{(g,h)}} \left( \frac{1}{e_{t} (\Nu (g,h))}\right).
\end{equation}
\end{Proposition}

\begin{proof}
Keeping in mind that $M^{(g,h)}$ is a compact manifold, the first
assertion follows by the argument originally due to
\cite{AtiyahSegal:IndexII}. 

Now examine the diagram
\[
\begin{CD}
 B\Lambda\times M^{(g,h)} @> i >> E\Lambda\times_{(g,h)} M @> j >>
  EG\times_{G} M \\
 @V 1 \times \pi^{M^{(g,h)}} VV @V \pi^{M,\Lambda} VV @VV \pi^{M,G}
V \\
B\Lambda @= B\Lambda @> B (g,h) >> BG
\end{CD}
\]
The right square is a pull-back, so
\[
     B (g,h)^{*}\pi_{t}^{M,G} = \pi_{t}^{M,\Lambda}j^{*}.
\]
It follows that
\[
    \hkr_{(g,h)} \pi_{!}^{M,G} (1) = \pi_{t}^{M,\Lambda} (1),
\]
considered as an element of $D$.

The fixed-point formula asserts that
\[
     \pi_{!}^{M,\Lambda} (1) = 1 \otimes \pi_{!}^{M^{(g,h)}} \left(
\frac{1}{e_{t} (\Nu (g,h))}\right)
\]
in $L$: but in fact we know that the left-hand side is an element of
$D\subset L$.  It follows that the right hand side is too, and
\[
    \hkr_{(g,h)} \pi_{!}^{M,G} (1) = 1 \otimes \pi_{!}^{M^{(g,h)}}
\left( \frac{1}{e_{t} (\Nu (g,h))}\right).
\]
The rest is easy.
\end{proof}

The formula \eqref{eq:4} is the analogue for the $t$-genus of the
orbifold elliptic genera of \cite{DMVV:egspsqs,BL:egsvoeg}.  To see
this, let $A$ be the (abelian) subgroup of $G$ generated by $g$ and
$h$, and suppose that $\Nu (g,h)$ decomposes as a sum
\[
   \Nu (g,h) \cong L_{1} \oplus \dotsb \oplus L_{r}
\]
of complex line bundles, with $A$ acting on $L_{i}$ by the character
$\chi_{i}$.  Let
\[
      e (\chi_{i}) \in E (BA)
\]
be the euler class of the character $\chi_{i}$, using the orientation
$t$, and let
\[
      y_{i} \in E (M^{(g,h)})
\]
be the (non-equivariant) euler class of the line bundle $L_{i}$.  Then
\[
   e_{t} (\Nu (g,h)) = \prod_{i} y_{i} \fs{F} e (\chi_{i}).
\]

We can be even more explicit.  Let
\[
   F (x,y) \in E \psb{x,y}
\]
be the formal group law over $E$ induced by the orientation $t$.  If
$R$ is a complete local $E$-algebra, let us write $F (R)$ for the
maximal ideal of $R$, considered as an abelian group using the power
series $F$ to perform addition.  

Associating to a character $\lambda\in \Lambda_{n}^{*}$ its first
chern class in $E$-theory using the orientation $t$ defines a
group homomorphism 
\[
   \Lambda_{n}^{*} \xra{} F (E (B\Lambda_{n})), 
\]
which gives rise to a homomorphism 
\[
 \Lambda_{n}^{*} \xra{} F (D_{n});
\]
in fact, this is the ``level structure'' referred to in
\S\ref{sec:character-theory}. 

The
dual of the epimorphism
\[
   \Lambda \rightarrow \Lambda_{n} \rightarrow A
\]
is a monomorphism
\[
    A^{*} \rightarrow \Lambda_{n}^{*}
\]
which composes with the level structure to give a homomorphism
\[
   \phi: A^{*} \rightarrow \Lambda_{n}^{*} \xra{} F (D_{n}).
\]
By construction,
\[
   \phi (\chi_{i}) = e (\chi_{i}).
\]
If $e_{1},e_{2}$ are a basis for $\Lambda_{n}^{*}$, then we can write
\[
    \chi_{i} = a_{i} e_{1} + b_{i} e_{2}
\]
in $\Lambda_{n}^{*}$, where $a_{i},b_{i}\in \Z/n$.  If $v_{i} = \ell
(e_{i})$ for $i=1,2$, then
\[
     \phi (\chi_{i}) = [a_{i}] (v_{1}) \fs{F} [b_{i}] (v_{2}).
\]
Our typical summand in the formula \eqref{eq:4} for the orbifold genus
becomes
\begin{equation} \label{eq:7}
1 \otimes \pi_{!}^{M^{(g,h)}} \left( \prod_{i} \frac{1}{y_{i}\fs{F}
[a_{i}] (v_{1}) \fs{F} [b_{i}] (v_{2})}\right).
\end{equation}

It is customary to calculate expressions like \eqref{eq:7} by using
the topological Riemann-Roch formula to pass to ordinary cohomology.
In fact this approach is not available in our situation. To do so, one
introduces the exponential
\[
    \exp : \Gah\to F
\]
of the group law $F$, and finds $x_{i}\in L\otimes E (M^{(g,h)})$ and
$w_{i}\in \Gah (L)$ such that
\begin{align*}
    y_{i} & = f (x_{i}) \\
    v_{1} & = f (w_{1}) \\
    v_{2} & = f (w_{2}).
\end{align*}
However, if $v_{1} = f (w_{1})$ then $0 = [n]_{F} (v_{1}) = f
(nw_{1})$, which implies that $nw_{1}=0$, and, as $L$ is torsion free,
we must have $w_{1} = 0$!

Nevertheless, we shall proceed formally in order to compare our
formula with those of \cite{DMVV:egspsqs,BL:egsvoeg}.  We have
\[
       e_{t} (\Nu (g,h)) = \prod_{i} f (x_{i} + a_{i} w_{1} + b_{i}
w_{2}).
\]
Let $u_{j}, j=1,\dotsc,r$ be the roots of the total Chern class of the
tangent bundle of $M^{(g,h)}$:
\[
   c (M^{(g,h)}) = \prod_{j} (1 - u_{j}).
\]
Then the Riemann-Roch formula gives
\begin{equation} \label{eq:8}
 \pi_{t}^{M^{(g,h)}} \left( \frac{1}{e_{t} (\Nu (g,h))}\right) =
\int_{M^{(g,h)}} \left( \prod_{j} \left(\frac{u_{j}}{f (u_{j})}
\right) \prod_{i}\frac{1}{f ( x_{i} + a_{i}w_{1} + b_{i} w_{2})}
\right).
\end{equation}

In \cite{BL:egsvoeg}, Borisov and Libgober use the two-variable
elliptic whose exponential is
\begin{equation} \label{eq:14}
        f (x) = \frac{\theta (x,\tau)}{\theta (x-z,\tau)},
\end{equation}
where
\begin{equation} \label{eq:10}
\theta (x,\tau) = -i q^{\tfrac{1}{8}} (e^{\frac{x}{2}} -
e^{-\frac{x}{2}})\prod_{n\geq 1} (1-q^{n}) (1-q^{n}e^{x})
(1-q^{n}e^{-x}),
\end{equation}
$\tau$ is a complex number with positive imaginary part, and
$q=e^{2\pi i \tau}$. (We have adopted slightly different conventions
regarding factors of $2\pi$.  The simplest way to compare is to say
that we work with the elliptic curve $\C/ (2\pi i \Z + 2\pi i \tau
\Z)$, while they work with the elliptic curve $\C/ (\Z+\tau \Z)$)

Their expression for the orbifold two-variable genus
is
\[
   \orbtvg (M,G) = \frac{1}{|G|} \sum_{gh=hg} \Phi_{g,h},
\]
where, with our conventions, 
\begin{equation} \label{eq:9}
   \Phi_{g,h} = \int_{M^{(g,h)}} 
\left( 
\prod_{j} \left(\frac{u_{j}}{f(u_{j})} 
\right) 
\prod_{i}\frac{1}
         {f ( x_{i} + A_{i} (1/n) - B_{i}(\tau/n))}
         e^{zB_{i}/n} 
\right).
\end{equation}
and the $A_{i}$ and $B_{i}$ are \emph{integer} representatives of
$a_{i}$ and $b_{i}$.

This differs from the formal expression \eqref{eq:8} by only the
factors $e^{zb_{i}/n}$.  These factors are familiar from the study of
equivariant genera; for example they are analogous to the factors
$\nu_{s}^{1/k}$ in (11.26) of \cite{BottTaubes:Rig}, $S (c_{1}
(\Nu^{1/n})_{\T})$ in (6.17) of \cite{AndoBasterra:WGEEC}, or
$u^{\frac{k}{n}\Ihat (\barm)}$ in (5.16) of \cite{Ando:AESO}.  Their
role is to make the expression \eqref{eq:9} independent of the choice
of representatives $A_{i}$ and $B_{i}$.  It is necessary to introduce
these factors because $f$ is not doubly periodic; instead we have
\begin{equation}  \label{eq:15}
  f (x+2\pi i \ell +  2\pi i k \tau)  = y^{-k} f (x),
\end{equation}
where $y=e^{z}$, as one checks easily using \eqref{eq:10}.  So the
expression doesn't depend on the choice of $A_{i}$.  If $B_{i}' =
B_{i} + n \delta_{i}$, then
\begin{align*}
   \prod_{i} f (x_{i} + A_{i}/n - B_{i}'\tau/n) e^{z B_{i}'/n}& =
\prod_{i} f (x_{i} + A_{i}/n - B_{i}\tau/n) y^{-\sum \delta_{i}}
   e^{z B_{i}/n} y^{\sum \delta_{i}} \\
& = \prod_{i} f (x_{i} + A_{i}/n - B_{i}\tau/n)e^{zB_{i}/n}.
\end{align*}
This is not an issue in our expression \eqref{eq:7}, and so the
factors $e^{zB_{i}/n}$ have no role in our genus.

\section{Comparison with the analytic equivariant genus}
\label{sec:comp-with-analyt}

In fact, we can use
the expression \eqref{eq:14} for the two-variable elliptic genus in
terms of theta functions to construct a Thom class in Grojnowski's equivariant
cohomology, following
\cite{Rosu:Rigidity,AndoBasterra:WGEEC,Ando:AESO}.  When $G=\T[n]$ is a
cyclic group of order $n$ acting on a compact manifold $M$, we can write down a
formula for a $\T$-equivariant genus on $\T\times_{G}M$, which 
by Shapiro's lemma is a sensible notion of $G$-equivariant genus on
$M$.  When we do so, we obtain the formula of \cite{DMVV:egspsqs,BL:egsvoeg}.

Let $\Lambda$ be the lattice $(2\pi i \Z + 2 \pi i \tau \Z)$, let
$C$ be the elliptic curve $\C/\Lambda$, and let $E_{\T}$ be
Grojnowski's equivariant elliptic cohomology associated to $C$.  
For convenience we identify 
\[
   \T\cong \R/\Z
\]
so that
\[
   G \cong \Z[\frac{1}{n}]/ \Z.
\]
We identify 
\[
    \T \times \T \cong C
\]
by the formula 
\begin{equation}
   (r+ \Z,s+\Z) \mapsto 2\pi i r + 2\pi i s \tau  + \Lambda.
\end{equation}

Let 
$M$ be a $G$-manifold with an equivariant complex structure on its
tangent bundle $T$.  We define 
\[
    E_{G} (M) \eqdef E_{\T} (\T\times_{G} M).
\]
For $a\in C$, 
\[
    (\T\times_{G}M)^{a} = 0
\]
unless $a\in C[n]$, so $E_{T} (\T\times_{G}M)$ is a metropolitan sheaf
(collection of skyscraper sheaves) supported at $C[n]$.  The stalk
at a point $a$ of order dividing $n$ is 
\[
     H (EG\times_{G} X^{a}).
\]

Let $\cochars$ be the lattice of cocharacters in $Spin (2d)$.  In
\cite{Ando:AESO}, the first author constructed orientations for theta
functions 
\[
\Theta = \Theta (\slot,\tau): \cochars\otimes \C \rightarrow \C
\]
satisfying 
\begin{equation} \label{eq:13}
   \Theta (x+2\pi i \ell+ 2 \pi  i  k \tau,\tau ) = 
   \exp (- 2\pi i I (k, x))
   \exp (-\phi (k))
   \Theta (x,\tau)
\end{equation}
for $x\in \cochars\otimes \C$ and $k, \ell\in \cochars$, 
where 
\begin{align*}
 \phi: & \cochars \to \Z \\
 I  : & \cochars \times\cochars \to \Z
\end{align*}
are respectively quadratic and bilinear functions related by 
\[
   \phi (\ell + \ell') = \phi (\ell) + I (\ell,\ell') + \phi (\ell').
\]

The building block of the orientation is a family of functions $F$
which we now describe.  
For simplicity we have supposed that $TM$ is a complex vector bundle,
and so our structure group is $U (d)$ instead of $Spin (2d)$.  We let
$T$ be the maximal torus of diagonal matrices; our choices so far
identify the lattice $\cochars=\hom (\T,T)$ of cocharacters with
$\Z^{d}$ in the usual way. 

Let $(g,h) \in G^{2}$; let $a$ be the corresponding point of $C[n]$,
and let $A\subseteq G$ be the subgroup 
generated by $g$ and $h$.  The action of $A$
on $TM\restr{M^{A}}$ is described by 
characters 
\[
     m = (m_{1},\dotsc ,m_{d}) \in \hom (A,T) \cong
     \cochars/|A|\cochars\cong (\Z/|A|)^{d}.
\]
Choose integer lifts 
\[
    \barm = (\barm_{1},\dotsc ,\barm_{d}) \in \cochars \cong \Z^{d}.
\]
Choose
\[
   \bara = 2\pi i \frac{\ell}{n}+ 2\pi i \frac{k}{n}\tau 
\]
so that 
\[
    a = \bara + \Lambda.
\]
In terms of these choices, the stalk of the orientation at
$(g,h)$ is built from the function 
\[
F (x)   = \exp (2\pi i \frac{k}{n}I (x,\barm))
                           \exp (2\pi i \frac{k}{n}\phi (\barm)\tau)
                           \Theta (x + \barm\otimes \bara)
\]
The essential feature of $F$ is that the functional equation
\eqref{eq:13} satisfied by $\Theta$ implies that $F$ is Weyl invariant
and independent of the choice of preimage $\barm$, and its
dependence on $\bara$ is under control.

Now consider the exponential $f$~\eqref{eq:14} associated to the two-variable
elliptic genus.  Comparison of its functional equation \eqref{eq:15}
with the functional equation \eqref{eq:13} for $\Theta$ suggests that we
should build the orientation for the two-variable genus from the
simpler function 
\begin{equation} \label{eq:16}
   F (x)  = 
   \exp (2\pi i \frac{k}{n}z\sum \barm_{j})
   \prod_{j} f (x_{j}+ 2\pi i \barm_{j}\ell/n + 
                       2\pi i \tau \barm_{j}k/n).
\end{equation}
The argument at the end of \S\ref{sec:orbifold-sigma-genus} shows
that indeed, this $F$ is independent of the choice of lift $\barm$.
In fact, the simple transformation rule \eqref{eq:15} for $f$ implies
that $F$ is also independent of the choice of representative $\bara$
for $a$.

Now use this $F$ to write down a class $\antvg (M,G) \in E_{G} (M)$
following the instructions in \cite{Ando:AESO}. (In general one gets a
section of the cohomology of the Thom space, but the Thom isomorphism
in ordinary cohomology defined by $f$ identifies this with
$E_{G} (M)$).  More precisely, the formula for the 
value in the stalk at $a\in C[n]$ is  the one for
$\gamma_{a}$ before Lemma 8.12, taking $V'$ to be trivial and
$\theta'$ to be $1$.  That formula refers to an expression $R$ which is
defined in terms of $F$ in Lemma 5.28. (The expression for $R$ also 
includes a  product of $\sigma$ functions, which should be replaced
with the corresponding product of $f$'s).

\begin{Theorem} 
With these substitutions, the value of $\antvg (M,G)$ in the stalk at
$(g,h)$ is the class in $H(EG \times_{G} M^{(g,h)})$ whose restriction to 
$H (BA\times M^{(g,h)})$ is the integrand in the
summand $\Phi_{g,-h}$ of the 
orbifold two-variable elliptic genus (see~\eqref{eq:9}). \qed
\end{Theorem}

With the remarks so far, the proof is straightforward.  We omit the
details, except note that if $(g,h) = (\ell/n+\Z,k/n+\Z)$, 
then in \eqref{eq:9},  $(A_{i},B_{i})$ can be taken to run over the set
$(\ell\barm_{i},k\barm_{i})$.  Thus 
the typical factor in \eqref{eq:9} can easily be seen to identify
with the typical factor in~\eqref{eq:16}.

The genus associated to $\antvg (M,G)$ is the global section of
$\O_{C[n]}$ whose value at $(g,h)$ is 
\[
    \int_{M^{(g,h)}} \antvg (M,G)_{g,h}.
\]
Summing over $(g,h)\in G$, we get exactly $|G|$ times the orbifold
two-variable elliptic genus of \cite{DMVV:egspsqs,BL:egsvoeg}.

\section{The cocycle} \label{sec:cocyle}

We now return to the orbifold sigma genus, and study the effect of
varying the $\busix$ structure.
The fibration of infinite loop spaces
\[
    K (\Z,3) \rightarrow \busix \rightarrow BSU
\]
gives a map of $\einfty$ ring spectra
\[
    \sinf K (\Z,3)_{+} \xra{i} \musix.
\]
If $\Ell{C} = (E,C,t)$ is an elliptic spectrum, then the sigma
orientation
\[
\sigma (\Ell{C}) : \musix \rightarrow E
\]
gives rise to a map of ring spectra
\[
   \sr{\Ell{C}} \eqdef \sigma (\Ell{C}) i: \sinf K (\Z,3)_{+} \xra{}
E.
\]
In particular, $\sr{\Ell{C}}$ is a Thom class for the trivial bundle
over $K (\Z,3)$, and so it is a \emph{unit} of $E^{0}K (\Z,3)$.

If $\alpha= (g,h): \Lambda\to G$, then we define
\begin{equation} \label{eq:12}
    \delta_{n} (\alpha) = \delta_{n} (u,\Ell{C},\alpha) \eqdef
\hkr_{g,h} (u^{*}\sr{\Ell{C}}) \in D^{\times}.
\end{equation}

\begin{Lemma} \label{t-le-delta}
If
\[
\alpha' = \Lambda_{n+1} \rightarrow \Lambda_{n} \xra{\alpha} G,
\]
then
\begin{equation} \label{eq:3}
     \delta_{n+1} (\alpha') = \delta_{n} (\alpha),
\end{equation}
and so we have a well-defined unit $\delta (g,h)\in D^{\times}$.  It
satisfies
\begin{align*}
    \delta (g,h)^{n} & = 1\\
    \delta (h,g) & = \delta (g,h)^{-1} \\
    \delta (g+g',h)  & = \delta (g,h) \delta (g',h) \\
    \delta (g,h+h') & = \delta_{n} (g,h) \delta (g',h)
\end{align*}
for any $g,g', h,$ and $h'$ in $G$ for which these equations make
sense.
\end{Lemma}

\begin{proof}
The arguments for the various claims are similar to each other; as an
illustration we show that $\delta$ is exponential in the first
variable, and for simplicity we suppose that $G$ is abelian.  The
proof in the general case will be given in \S\ref{sec-non-abelian}.

Let $C$ be the cyclic group of order $n$.  The universal example of an
abelian group with three elements $(g_{1}, g_{2}, h)$ of order $n$ is
$C^{3}$.  Let 
\begin{align*}
    C^{2} & \xra{(g_{1},h)} C^{3}\\
    C^{2} & \xra{(g_{2},h)} C^{3}\\
    C^{2} & \xra{(g_{1}+g_{2},h)} C^{3}
\end{align*}
be the maps which represent the selection of the indicated pairs elements
formed from the triple $(g_{1},g_{2},h)$.  It suffices to show that for
every homotopy class
\[
     u: BC^{3} \xra{} K (\Z,3),
\]
the outside rectangle of the diagram
\begin{equation} \label{eq:20}
\begin{CD}
BC^{2} @>(g_{1}+g_{2},h)>> BC^{3} @>u>> K (\Z,3) @>\sigma(\Ell{C})>> E\\
@V \Delta VV               @.              @AAA                    @AAA \\
BC^{2}\times BC^{2} @> (g_{1},h)\times (g_{2},h) >> BC^{3} \times
BC^{3} @> u \times u >> K (\Z,3)\times K (\Z,3) @>\sigma
(\Ell{C})\times\sigma (\Ell{C}) >> E\times E
\end{CD}
\end{equation}
commutes.  The right-side rectangle commutes, because
$\sigma(\Ell{C})$ is a map of ring spectra.  To show that the
left-side rectangle commutes, consider the diagram
\begin{equation} \label{eq:19}
\begin{CD}
BC^{2} @>(g_{1}+g_{2},h)>> BC^{3} @>v>> K (C,2) @> \beta >> K (\Z,3)\\
@V \Delta VV               @.              @AAA                    @AAA \\
BC^{2}\times BC^{2} @> (g_{1},h)\times (g_{2},h) >> BC^{3} \times
BC^{3} @> v \times v >> K (C,2)\times K (C,2) @> \beta\times \beta >>
K (\Z,3)\times K (\Z,3).
\end{CD}
\end{equation}
Once again, the right square commutes, this time because the Bockstein
is an additive group homomorphism.

An easy calculation shows that the Bockstein
\[
    \beta: H^{2} (BC^{3};C) \rightarrow H^{3} (BC^{3};\Z)
\]
is surjective, so there is a $v$ such that $\beta v = u$.  Moreover
\[
    H^{2} (BC^{3};C)\cong C\{\mu_{12},\mu_{13},\mu_{23}\},
\]
where
\[
\mu_{ij}: BC^{3}\cong K (C,1)^{3}\rightarrow K (C,2)
\]
is the map which represents the natural transformation
\[
     \mu_{ij}: H^{1} (X;C)^{3} \rightarrow H^{2} (C;C)
\]
given by
\[
    \mu_{ij} (x_{1},x_{2},x_{3}) = x_{i}\cup x_{j}.
\]
If
\[
    v = a \mu_{12} + b \mu_{13} + c \mu_{23},
\]
then the top row of the diagram represents the natural transformation
\[
   H^{1} (X;C)^{2} \xra{} H^{3} (X;\Z)
\]
given by
\[
    (x,y) \mapsto \beta ( a x \cup x + b x \cup y + c x\cup y),
\]
while the other (i.e. counterclockwise) composition represents the
natural transformation 
\[
    (x,y) \mapsto \beta ( b x \cup y + c x \cup y).
\]
These coincide since $\beta (x\cup x) = 0.$

In other words, the diagram \eqref{eq:19} commutes for \emph{any}
$v$, and so the diagram \eqref{eq:20} does too, as required.
\end{proof}

\subsection{The Weil pairing}

As an example, let's consider the case that $G=\Lambda$, and
\[
    u: BG \rightarrow K (\Z/N,2)
\]
is the map representing the cup product: indeed $u$ is a generator of
$E^{0}BG$.  A homomorphism
\[
           \alpha: \Lambda \to G
\]
gives rise to a homomorphism
\[
           G^{*} \xra{} \Lambda^{*} \rightarrow C[N].
\]
Thus we may view the homomorphism $\alpha= (g,h)$ as a pair of
$N$-torsion points of $C[N]$.  The argument of \cite{AS:WPMK} then
shows

\begin{Lemma}
$\delta (u)$ is the \emph{Weil pairing} of the elliptic curve
$C$. \qed
\end{Lemma}


\section{Discrete torsion} \label{sec:formula}

We are now ready to state our basic formula.

\begin{Theorem}
If $\alpha = (g,h): \Lambda_{n} \rightarrow G$, then
\[
         \hkr_{\alpha} \sigmagenus (M,\ell+\iota u \pi,\Ell{C})_{G} =
\delta (u,\Ell{C},\alpha) \hkr_{\alpha} \sigmagenus
(M,\ell,\Ell{C})_{G},
\]
and so abbreviating $\delta (g,h) = \delta (u,\Ell{C},\alpha)$, we
have
\[
        \orbsigma (M,\ell+\iota u \pi,\Ell{C})_{G} = \sum_{gh=hg}
\delta (g,h) \hkr_{g,h} \sigmagenus (M,\ell,\Ell{C})_{G} \in \pi_{-d}
E.
\]
\end{Theorem}

\begin{proof}
Since
\[
   \sigma (\Ell{C}): \musix \rightarrow E
\]
is a map of ring spectra, and since the multiplication on $\musix$
arises from the addition on $\busix$, we have
\begin{align*}
      \SigmaClass  (M,\ell+\iota u\pi,\Ell{C})_{G} & =
\SigmaClass  (M,\ell,\Ell{C})_{G} \SigmaClass (M,\iota u\pi,\Ell{C})_{G}  \\
& = \SigmaClass (M,\ell,\Ell{C})_{G} \pi^{*} u^{*}\sr{\Ell{C}}.
\end{align*}

Recall that $E^{0} (\borel{V}-\borel{T})$ is an $E^{0}
(\borel{M})$-module, and so an $E^{0} (BG)$-module via $E^{0} (\pi)$.
As such the Pontrjagin-Thom map
\[
     \pontthom{V}_{G}: E^{0} (\borel{V} - \borel{T}) \rightarrow E^{0} (BG)
\]
is a homomorphism of $E^{0} (BG)$-modules. It follows that
\begin{align*}
     \sigmagenus (M,\ell+\iota u \pi,\Ell{C})_{G} & = 
     \pontthom{V}_{G}(\SigmaClass  (M,\ell+\iota u\pi,\Ell{C})_{G}) \\
    & = \pontthom{V}_{G} 
       (\SigmaClass  (M,\ell,\Ell{C})_{G} \pi^{*} u^{*}\sr{\Ell{C}}) \\
     & = u^{*}\sr{\Ell{C}} \sigmagenus  (M,\ell,\Ell{C})_{G}
\end{align*}
in $E^{*} (BG)$.  If
\[
\alpha : \Lambda_{n} \to G,
\]
then applying $\hkr_{\alpha}$ gives
\begin{align*}
 \hkr_{\alpha} \sigmagenus (M,\ell+\iota u\pi,\Ell{C})_{G} & =
     \hkr_{\alpha} u^{*}\sr{\Ell{C}} 
     \hkr_{\alpha}\sigmagenus  (M,\ell,\Ell{C})_{G} \\  
     & = \delta_{n} (u,\Ell{C},\alpha) 
     \hkr_{\alpha}\sigmagenus (M,\ell,\Ell{C})_{G}.
\end{align*}
as required.
\end{proof}

\section{The non-abelian Case}\label{sec-non-abelian}

In this section, we prove Lemma \ref{t-le-delta} in the case that $G$
is non-abelian.  We fix an $n$ sufficiently large that $|G|$ divides
$n$. We first construct an isomorphism $H^3(B\Lambda_n)\cong \Z/n.$

Consider the following commutative diagram, where the columns are
universal coefficient exact sequences.

$$\xymatrix{{\Ext^1(H_1(B\Lambda_n);\Z)}\ar[r]\ar[d] &
{\Ext^1(H_1(B\Lambda_n);\Z/n)}\ar[d] \\ {H^2(B\Lambda_n)}\ar[r]\ar[d]
&
{H^2(B\Lambda_n;\Z/n)}\ar[r]\ar[d] & 
{H^3(B\Lambda_n)}\ar@{-->}[dl]^{\cong} \\
{\Hom(H_2(B\Lambda_n);\Z)}\ar[r] & {\Hom(H_2(B\Lambda_n);\Z/n)}.}$$

Here, the middle row is part of Bockstein long exact sequence.  Note
that this is a short exact sequence since multiplication by $n$ kills
group cohomology of $\Lambda_n$. Since $H_2(B\Lambda_n)$ is torsion,
$\Hom(H_2(B\Lambda_n);\Z)=0$.  We therefore obtain the dotted arrow.
Since $\Ext^1(H_1(B\Lambda_n);-)$ is right exact,the top map is a
surjection.  It follows that the dotted arrow is an isomorphism.  Now
it is not hard to check that $e_1\otimes e_2-e_2\otimes e_1$ is a
generator for $H_2(B\Lambda_n)\cong \Z/n$.  Composing the dotted arrow
above with evaluation on this generator yields an
isomorphism $$H^3(B\Lambda_n)\cong \Z/n.$$

\begin{Definition} \label{def-epsilon} If $u$ is an element in
$H^3(BG)$ and $\alpha=(g,h):\Lambda_n\to G$ is any map, then we obtain
a class $\alpha^*u\in H^3(B\Lambda_n)$.  Let $\epsilon (g,h) =
\epsilon_u^n(g,h)$ be the image of this class in $\Z/n$ under the
isomorphism above.
\end{Definition}

\begin{Remark} It is easy to check that
$$\epsilon(g,h)=\tilde u(g,h)-\tilde u(h,g)$$ where $\tilde u:G\times
G\to \Z/n$ is a 2-cocycle whose cohomology class maps to $u$ under the
Bockstein.
\end{Remark}

\begin{Lemma}\label{lem:cocycle} Whenever the expressions are defined,
the following properties hold.
$$\epsilon(g,h)=-\epsilon(h,g)$$
$$\epsilon(h,j)-\epsilon(gh,j)+\epsilon(g,hj)-\epsilon(g,h)=0$$
$$\epsilon(gg',h)=\epsilon(g,h)+\epsilon(g',h)$$
$$\epsilon(g,hh')=\epsilon(g,h)+\epsilon(g,h')$$
\end{Lemma}

\begin{proof} The first and second properties follow easily from the
remark.  For the third property, it suffices to show that
\[
\tilde u (gg',h)-\tilde u (h,gg')-\tilde u (g,h)+\tilde u (h,g)-\tilde
u (g',h)+\tilde u (h,g')
\] is zero.  Since $\tilde u$ is a cocycle, we may rewrite the
expression using the following equations:
\begin{gather*}\tilde u
(gg',h)-\tilde u (g',h)=\tilde u (g,g'h)-\tilde u (g,g'),\\
-\tilde u (h,gg')+\tilde u (h,g)= \tilde u (g,g')-\tilde u
(hg,g'),\\
-\tilde u (g,h)+\tilde u (h,g')=\tilde u
(gh,g')-\tilde u (g,hg').
\end{gather*}
Then, canceling terms, we get
\[
\tilde u (g,g'h)+\tilde u (gh,g')-\tilde u (hg,g')-\tilde u (g,hg').
\]
Since $gh=hg$ and $g'h=hg'$, this is zero as needed.

The last property follows similarly, or directly from the first and
third.
\end{proof}

\begin{proof}[Proof of Lemma \ref{t-le-delta}] Let
\[
F: B\Lambda_{n} \rightarrow K (\Z,3)
\]
represent the element $1\in \Z/n\cong H^{3} (B\Lambda_{n})$ under the
isomorphism above.  Let $x\in D$ be the image of
\[
    F^{*} \sr{\Ell{C}} \in E^{0}B\Lambda_{n}
\]
under the map tautological map 
\[
    E^{0}B\Lambda_{n} \xra{}D_{n} \xra{} D.
\]
It is easy to check from the definitions~\eqref{eq:12}
and~\eqref{def-epsilon} of $\delta$ and $\epsilon$ that
\[
   \delta(g,h) = x^{\epsilon_{g,h}}.
\]
The result now follows from Lemma \ref{lem:cocycle}.
\end{proof}


\begin{thebibliography}{DMVV97}

\bibitem[AB02]{AndoBasterra:WGEEC}
Matthew Ando and Maria Basterra.
\newblock The {W}itten genus and equivariant elliptic cohomology.
\newblock {\em Mathematische Zeitschrift}, 240(4):787--822, 2002.

\bibitem[AHS01]{AHS:ESWGTC}
Matthew Ando, Michael~J. Hopkins, and Neil~P. Strickland.
\newblock Elliptic spectra, the {W}itten genus, and the theorem of the cube.
\newblock {\em Inventiones Mathematicae}, 146:595--687, 2001.
\newblock DOI 10.1007/s002220100175.

\bibitem[AHS03]{AHS:hinfty}
Matthew Ando, Michael~J. Hopkins, and Neil~P. Strickland.
\newblock The sigma orientation is an ${H}_\infty$ map.
\newblock {\em Amer. J. Math.}, To appear.

\bibitem[And03]{Ando:AESO}
Matthew Ando.
\newblock The sigma orientation for analytic circle-equivariant elliptic
  cohomology.
\newblock {\em Geometry and Topology}, 7:91--153, 2003.
\newblock math.AT/0201092.

\bibitem[AS68]{AtiyahSegal:IndexII}
M.~F. Atiyah and G.~B. Segal.
\newblock The index of elliptic operators. {I}{I}.
\newblock {\em Ann. of Math. (2)}, 87:531--545, 1968.

\bibitem[AS01]{AS:WPMK}
M.~Ando and N.~P. Strickland.
\newblock Weil pairings and {M}orava {$K$}-theory.
\newblock {\em Topology}, 40(1):127--156, 2001.

\bibitem[BL02]{BL:egsvoeg}
Lev~A. Borisov and Anatoly Libgober.
\newblock Elliptic genera of singular varieties, orbifold elliptic genus and
  chiral de {R}ham complex.
\newblock In {\em Mirror symmetry, IV (Montreal, QC, 2000)}, volume~33 of {\em
  AMS/IP Stud. Adv. Math.}, pages 325--342. Amer. Math. Soc., Providence, RI,
  2002.

\bibitem[BM94]{bm:gd4ccI}
J.-L. Brylinski and D.~A. McLaughlin.
\newblock The geometry of degree-four characteristic classes and of line
  bundles on loop spaces. {I}.
\newblock {\em Duke Math. J.}, 75(3):603--638, 1994.

\bibitem[BT89]{BottTaubes:Rig}
Raoul Bott and Clifford Taubes.
\newblock On the rigidity theorems of {W}itten.
\newblock {\em J. of the Amer. Math. Soc.}, 2, 1989.

\bibitem[DMVV97]{DMVV:egspsqs}
Robbert Dijkgraaf, Gregory Moore, Erik Verlinde, and Herman Verlinde.
\newblock Elliptic genera of symmetric products and second quantized strings.
\newblock {\em Comm. Math. Phys.}, 185(1):197--209, 1997.

\bibitem[Dri74]{Drinfeld:EM}
V.~G. Drinfeld.
\newblock Elliptic modules.
\newblock {\em Math. USSR-Sb.}, 23(4):561--592, 1974.

\bibitem[EOTY89]{EOTY}
H.~Eguchi, H.~Ooguri, A.~Taormina, and S.-K. Yang.
\newblock Superconformal algebras and string compactification on manifolds with
  ${S}{U}({N})$ holonomy.
\newblock {\em Nucl. Phys. B}, 315, 1989.

\bibitem[GMS00]{GMS:Gerbes}
Vassily Gorbounov, Fyodor Malikov, and Vadim Schechtman.
\newblock Gerbes of chiral differential operators.
\newblock {\em Math. Res. Lett.}, 7(1):55--66, 2000.
\newblock math.AG/9906117.

\bibitem[HKR00]{HKR:ggc}
Michael~J. Hopkins, Nicholas~J. Kuhn, and Douglas~C. Ravenel.
\newblock Generalized group characters and complex oriented cohomology
  theories.
\newblock {\em J. Amer. Math. Soc.}, 13(3):553--594 (electronic), 2000.

\bibitem[Hop95]{ho:icm}
Michael~J. Hopkins.
\newblock Topological modular forms, the {W}itten genus, and the theorem of the
  cube.
\newblock In {\em Proceedings of the International Congress of Mathematicians,
  Vol.\ 1, 2 (Z\"urich, 1994)}, pages 554--565, Basel, 1995. Birkh\"auser.

\bibitem[Kri90]{MR91e:57059}
Igor~M. Krichever.
\newblock Generalized elliptic genera and {B}aker-{A}khiezer functions.
\newblock {\em Mat. Zametki}, 47(2):34--45, 158, 1990.

\bibitem[Ros01]{Rosu:Rigidity}
Ioanid Rosu.
\newblock Equivariant elliptic cohomology and rigidity.
\newblock {\em Amer. J. Math.}, 123(4):647--677, 2001.

\bibitem[Sha]{Sharpe:RDDT}
Eric Sharpe.
\newblock Recent developments in discrete torsion.
\newblock hep-th/0008191.

\bibitem[Str97]{Strickland:FiniteSubgps}
Neil~P. Strickland.
\newblock Finite subgroups of formal groups.
\newblock {\em J. Pure and Applied Algebra}, 121:161--208, 1997.

\bibitem[Vaf85]{vafa:dt}
C.~Vafa.
\newblock Modular invariance and discrete torsion.
\newblock {\em Nuclear Physics B}, 261:678--686, 1985.

\bibitem[Wit87]{Witten:EllQFT}
Edward Witten.
\newblock Elliptic genera and quantum field theory.
\newblock {\em Comm. Math. Phys.}, 109, 1987.

\end{thebibliography}

\end{document}